\documentclass[10pt]{amsart}
\usepackage{amssymb}
\begin{document}
\textwidth=5.9truein
\leftmargin=0.3truein
\rightmargin=0.3truein
\textheight=8.3truein
\title
{An Application of Expanders to 
${\mathbb B}(\ell_2)\otimes{\mathbb B}(\ell_2)$} 

\author{Narutaka OZAWA}

\address{Department of Mathematical Science, 
University of Tokyo, 153-8914, Japan}

\email{narutaka@ms.u-tokyo.ac.jp} 

\thanks{Partially supported by JSPS Postdoctral 
Fellowships for Research Abroad}

\begin{abstract}
With the help of Kirchberg's and Selberg's theorems, 
we prove that the minimal tensor product of ${\mathbb B}(\ell_2)$ 
with itself does not have the WEP (weak expectation property) of Lance. 
\end{abstract}

\newtheorem{thm}{Theorem}
\newtheorem{prop}[thm]{Proposition} 
\newtheorem{lem}[thm]{Lemma} 
\newtheorem{thm2}{Theorem}[section]
\newtheorem{cor2}[thm2]{Corollary} 
\newtheorem{lem2}[thm2]{Lemma} 
\theoremstyle{definition} 
\newtheorem{defn}[thm]{Definition} 
\newcommand{\N}{{\mathbb N}} 
\newcommand{\Z}{{\mathbb Z}} 
\newcommand{\C}{{\mathbb C}} 
\newcommand{\R}{{\mathbb R}} 
\newcommand{\F}{{\mathbb F}} 
\newcommand{\G}{\Gamma} 
\newcommand{\Lp}{\Lambda_p} 
\newcommand{\Lq}{\Lambda_q} 
\newcommand{\Lo}{\Lambda_\omega} 
\newcommand{\M}{{\mathbb M}}
\newcommand{\B}{{\mathbb B}} 
\newcommand{\hh}{{\mathcal H}} 
\newcommand{\hk}{{\mathcal K}} 
\newcommand{\id}{\mathrm{id}} 
\newcommand{\Ad}{\mathop{\mathrm{Ad}}} 
\newcommand{\Tr}{\mathop{\mathrm{Tr}}} 
\newcommand{\Sg}{\mathcal{S}} 
\newcommand{\OS}{OS} 
\newcommand{\e}{\varepsilon} 
\newcommand{\p}{\varphi} 
\newcommand{\spec}{\mathop{\mathrm{Sp}}} 
\newcommand{\ox}{\mathop{\otimes}}
\maketitle
\section{Introduction}
Kirchberg \cite{kirchberg} showed a remarkable theorem 
that there is a unique $C^*$-norm on the tensor product 
between a $C^*$-algebra with the LLP 
and a $C^*$-algebra with the WEP. 
(See \cite{pisier} for a simpler proof.) 
In the same paper, he raised several interesting problems. 
Among others, he asked if there is a unique $C^*$-norm 
on the tensor product of $\B(\ell_2)$ with itself. 
This problem was solved negatively by Junge and Pisier \cite{jp}. 
Their second approach uses the expanders (see also \cite{valette}). 
We refer the reader to a book of Lubotzky \cite{lubotzky} 
for the information of expanders. 
In this paper, we will give another application of expanders 
(or more precisely, of Selberg's theorem \cite{selberg}) 
to the tensor product of $\B(\ell_2)$ with itself. 
Our proof proceeds in the same spirit as that of Junge and 
Pisier \cite{jp} (and also of Voiculescu \cite{voiculescu}) 
to produce an uncountable family of operator spaces inside 
a separable metric space of operator spaces embeddable into 
the full group $C^*$-algebra $C^*(\F)$ of a free group $\F$. 
See \cite{jp} for the detail. 
\begin{thm}\label{thm} 
The $C^*$-algebra $\B(\ell_2)\ox_{\min}\B(\ell_2)$ 
does not have the WEP. 
\end{thm}
This theorem is a corollary of the following proposition, 
which is of independent interest.
\begin{prop}\label{prop} 
There are a set $\Lambda$ and an action $\sigma$ of $\G=PSL(2,\Z)$ 
on $\Lambda$ (as permutations) such that 
the corresponding full crossed product
$\ell_\infty(\Lambda)\rtimes\G$ 
does not have the LLP. 
\end{prop} 
We recall that an action $\alpha$ of a discrete group $\G$ 
on a $C^*$-algebra $A$ is a homomorphism of $\G$ into 
the group of $*$-automorphisms of $A$. 
A $C^*$-algebra with an action of $\G$ is called a $\G$-$C^*$-algebra 
and a map between $\G$-$C^*$-algebras is said to be $\G$-equivariant 
if the map is compatible with the $\G$-actions. 
The full crossed product $C^*$-algebra $A\rtimes_\alpha\G$ is 
then defined as the universal $C^*$-algebra generated by 
a copy of $A$ and a unitary representation $U$ of $\G$ 
under the relation $\Ad U(g)(a)=U(g)aU(g)^*=\alpha(g)(a)$ 
for $g\in\G$ and $a\in A$. 
We will often omit $\alpha$ and denote 
the full crossed $C^*$-algebra simply by $A\rtimes\G$. 
It follows from the $\G$-equivariant Stinespring theorem that 
a $\G$-equivariant unital completely positive map between 
$\G$-$C^*$-algebras naturally extends to a unital completely 
positive map between their full crossed products. 
This fact will be used in the proof of Lemma \ref{rep}. 

We recall the definitions of the LLP \cite{kirchberg} 
and the WEP \cite{lance}. 
\begin{defn} 
A unital $C^*$-algebra $A$ has the LLP (the local lifting property) 
if for any unital completely positive map 
$\phi$ from $A$ into a quotient $C^*$-algebra $B/J$ 
and any finite dimensional operator subsystem $E$ in $A$, 
there is a unital completely positive lifting of $\psi\colon E\to B$ 
of $\phi|_E$. 
A $C^*$-algebra $A$ has the WEP (the weak expectation property) 
if for any faithful representation $A\subset\B(\hh)$, 
there is a unital completely positive map 
$\Phi$ from $\B(\hh)$ into $A^{**}$ which is identical on $A$. 
\end{defn} 

Besides nuclear $C^*$-algebras, 
a typical example of a $C^*$-algebra with the LLP is 
the full group $C^*$-algebra $C^*(\F_\infty)$ of the 
free group $\F_\infty$ on countably many generators and 
a typical example with the WEP is $\B(\ell_2)$, 
the $C^*$-algebra of all bounded linear operators on the separable 
infinite dimensional Hilbert space $\ell_2$. 
Pisier \cite{pisier} showed that the LLP is closed 
under taking a full free product (see also \cite{boca}). 
It follows that the full group $C^*$-algebra $C^*(\G)$ has the LLP 
since $\G$ is isomorphic to the free product $(\Z/2\Z)*(\Z/3\Z)$. 
As we mentioned in the beginning, 
Kirchberg \cite{kirchberg} proved the following. 
\begin{thm}[Kirchberg \cite{kirchberg}]\label{kirchberg}
For $C^*$-algebras $A$ and $B$, we have 
\begin{enumerate}
\item $A\ox_{\min}B=A\ox_{\max}B$ 
if $A$ has the LLP and $B$ has the WEP, 
\item $A\ox_{\min}\B(\ell_2)=A\ox_{\max}\B(\ell_2)$ 
if and only if $A$ has the LLP, 
\item $C^*(\F_\infty)\ox_{\min}B=C^*(\F_\infty)\ox_{\max}B$ 
if and only if $B$ has the WEP. 
\end{enumerate}
\end{thm}

We use the following variant of 
the deep theorem of Selberg \cite{selberg} 
which has already been applied to $C^*$-algebras by Bekka \cite{bekka} 
to show that some full group $C^*$-algebras of residually finite 
groups are not residually finite dimensional $C^*$-algebras. 
We refer the reader to Sections 4.3 and 4.4 in \cite{lubotzky} 
for this theorem. 
\begin{thm}[Selberg \cite{selberg}]\label{selberg} 
The trivial representation of $SL(2,\Z)$ is isolated 
in the set of all unitary representations which factor 
through $SL(2,\Z/m\Z)$ for some $m\in\N$. 
\end{thm}
This theorem means that for any finite generating subset 
$\Sg$ of $SL(2,\Z)$, 
there are a constant $\kappa_{\Sg}>0$ and a continuous function 
$\alpha_{\Sg}\colon \R_{\geq0}\to\R_{\geq0}$ with $\alpha_{\Sg}(0)=0$ 
such that the following holds; 
if $\pi$ is a unitary representation of $SL(2,\Z)$ 
on a Hilbert space $\hh$, which factors through $SL(2,\Z/m\Z)$ 
for some $m\in\N$, and $\xi\in\hh$ is a unit vector with 
$\e=\max_{g\in\Sg}\|\pi_g\xi-\xi\|<\kappa_{\Sg}$, 
then there is a unit vector $\eta\in\hh$ 
such that $\pi_g\eta=\eta$ for all $g\in SL(2,\Z)$ 
and $\|\xi-\eta\|<\alpha_{\Sg}(\e)$. 
We observe here that the uniform convexity of a Hilbert space implies 
that, for each $n$, there is a continuous function 
$\beta_n\colon\R_{\geq0}\to\R_{\geq0}$ with $\beta_n(0)=0$, such that 
the following holds: 
if $\xi_1,\ldots,\xi_n$ are vectors in $\hh$
with $\|\xi_i\|\le 1$ such that $\|\sum_{i=1}^n\xi_i\|> n(1-\e)$, 
then $\|\xi_i-\xi_1\|<\beta_n(\e)$ for all $i$.  
In particular, if $u_1=1,u_2,\ldots,u_n$ are contractions on $\hh$
with $\|\sum_{i=1}^n u_i\|> n(1-\e)$, 
then there is a unit vector $\xi$ such that 
$\|u_i\xi-\xi\|<\beta_n(\e)$ for all $i$. 

\noindent{\bf Acknowledgment.}
The author thanks Professor Eberhard Kirchberg 
for pointing out a mistake in an earlier draft, 
Professor Gilles Pisier for suggesting a clearer presentation 
and the referee for some simplification of the proof of 
Theorem \ref{kex}. 
A part of this study was done under the support of JSPS Postdoctoral 
Fellowships for Research Abroad. 
\section{Proofs}\label{sec2} 
We recall from \cite{pisierex}\cite{jp} that 
$\OS_d$ is the set of all $d$-dimensional operator spaces, 
equipped with the cb Banach Mazur distance topology. 
By the definition of the LLP, the set of $d$-dimensional 
operator subspaces of a (not necessarily separable) 
$C^*$-algebra with the LLP is contained in 
the set of $d$-dimensional operator subspaces of 
the separable $C^*$-algebra $C^*(\F_\infty)$, and a fortiori, 
is separable in $\OS_d$. 
Therefore, to show a $C^*$-algebra $A$ does not have the LLP, 
it suffices to show that the set of $d$-dimensional operator 
subspaces of $A$ is not separable for some $d$. 
This was done by Junge and Pisier \cite{jp} for $A=\B(\ell_2)$. 
In this section, we will find an explicit example of 
an action $\sigma$ of $\G=PSL(2,\Z)$ on a set $\Lambda$ and 
$d$ such that the set of $d$-dimensional operator subspaces of 
$\ell_\infty(\Lambda)\rtimes\G$ is non-separable. 

For each prime number $p$, we let 
$\Lp$ be the projective space $(\Z_p^2-\{0\})/\Z_p^\times$, 
where $\Z_p^2$ is the two dimensional vector space over 
the finite field $\Z_p=\Z/p\Z$. 
We observe $|\Lp|=(p^2-1)/(p-1)=p+1$ and denote by $[t]$ 
the equivalence class of $(t\ 1)^T$ 
for $t\in\Z_p$ and by $[\infty]$ the equivalence class of $(1\ 0)^T$. 
The action of $\G$ on $\Z_p^2$ 
(through linear transformation by $SL(2,\Z_p)$) 
induces a transitive action $\sigma_p$ of $\G$ 
on the set $\Lp$. 
Let $\pi_p$ be the corresponding unitary representation 
of $\G$ on $\ell_2(\Lp)$. 
For example, letting 
$$h=\left(\begin{array}{ccc}1 & 1\\ 0 & 1\end{array}\right)\quad
\mbox{ and }\quad 
k=\left(\begin{array}{ccc}0 & 1\\ -1 & 0\end{array}\right)$$ 
be elements in $\G$, 
we have $\pi_p(h)\delta_{[t]}=\delta_{[t+1]}$ 
and $\pi_p(k)\delta_{[t]}=\delta_{[-t^{-1}]}$ 
for all $t$ in $\Z_p\cup\infty$. 
Let $\chi_p=|\Lp|^{-1/2}\sum_{x\in\Lp}\delta_x\in\ell_2(\Lp)$ 
be the constant function of norm $1$ and let $z\in\B(\ell_2(\Lp))$ 
be the projection onto the one dimensional subspace $\C\chi_p$. 
We observe that $z$ is in the center of $C^*(\pi_p(\G))$. 
Let $\pi_p^\circ$ be the subrepresentation of $\pi_p$, 
which is the restriction to the subspace $(1-z)\ell_2(\Lp)$. 
\begin{lem}\label{irr}
The representation $\pi_p^\circ$ is irreducible 
for every prime number $p$. 
\end{lem} 
\begin{proof}
This is well-known (cf.\ p71--72 in \cite{fh}), 
but we include the proof for the reader's convenience. 
We assert that the eigenspace $F$ of $\pi_p(h)$ 
w.r.t.\ the eigenvalue 1 is 
the two dimensional subspace $F'$ 
spanned by $\chi_p$ and $\delta_{[\infty]}$. 
Indeed, if 
$\zeta=\sum_{t\in\Z_p\cup\infty}c(t)\delta_{[t]}$ 
is in $F$, then it follows from the equation $\pi_p(h)\zeta=\zeta$ 
that $c(t)=c(0)$ for all $t\in\Z_p$. 
This shows $F\subset F'$, while the converse inclusion is clear. 
Therefore, $C^*(\pi_p^\circ(\G))$ contains 
the rank one projection onto the subspace spanned by 
$(1-z)\delta_{[\infty]}$ 
and the irreducibility of $\pi_p^\circ$ follows 
from the transitivity of $\sigma_p$. 
\end{proof}

Combined with Schur's lemma, this implies that, for $p\neq q$, 
any fixed vector for the representation $\overline{\pi_q}\otimes\pi_p$ 
is a scalar multiple of $\overline{\chi_q}\otimes\chi_p$, 
where $\overline{\pi_q}$ is the conjugate representation of 
$\pi_q$ on the conjugate Hilbert space $\overline{\ell_2(\Lq)}$. 

Let $\Omega$ be the set of all prime numbers and 
let $\Lambda=\bigsqcup_{p\in\Omega}\Lp$ be the disjoint union. 
Then, the collection $(\sigma_p)_{p\in\Omega}$ induces 
an action $\sigma$ of $\G$ on the set $\Lambda$ and 
an action $\alpha$ on $\ell_\infty(\Lambda)$. 
We denote by $U(g)$ the implementing unitary of $g\in\G$ 
in the full crossed product $\ell_\infty(\Lambda)\rtimes_\alpha\G$. 
Fixing a faithful representation $C^*(\G)\subset\B(\hh)$ of the 
full group $C^*$-algebra of $\G$, with $u(g)$ denoting the unitary 
corresponding to $g\in\G$, we define a covariant representation 
$$\rho\colon\ell_\infty(\Lambda)\rtimes\G
\to\B(\ell_2(\Lambda)\otimes\hh)$$ 
by 
$\ell_\infty(\Lambda)\ni a
\mapsto a\otimes 1\in\B(\ell_2(\Lambda)\otimes\hh)$ 
and 
$\G\ni g\mapsto\pi(g)\otimes u(g)\in\B(\ell_2(\Lambda)\otimes\hh)$, 
where we put $\pi(g)=\bigoplus_{p\in\Omega}\pi_p(g)\in\B(\ell_2(\Lambda))$. 

\begin{lem}\label{rep}
The representation $\rho$ is faithful. 
\end{lem}
\begin{proof}
Consider the diagonal embedding 
of $(\ell_\infty(\Lambda),\alpha)$ into 
$(\B(\ell_2(\Lambda)),\Ad \pi)$. 
Since $\ell_\infty(\Lambda)$ is the range of a $\G$-equivariant 
conditional expectation from $\B(\ell_2(\Lambda)$, 
the canonical morphism 
$$\ell_\infty(\Lambda)\mathop{\rtimes}_\alpha\G\to 
\B(\ell_2(\Lambda))\mathop{\rtimes}_{\Ad \pi}\G$$ 
is faithful. 
Since $\Ad \pi$ is inner, we have an isomorphism 
$$\B(\ell_2(\Lambda))\mathop{\rtimes}_{\Ad \pi}\G
\cong\B(\ell_2(\Lambda))\ox_{\max}C^*(\G),$$ 
where the implementing unitaries $U(g)$ in the left hand side 
are mapped to $\pi(g)\otimes u(g)$ in the right hand side. 
It follows from the LLP of $C^*(\G)$ and Theorem \ref{kirchberg} that 
$$\B(\ell_2(\Lambda))\ox_{\max}C^*(\G) 
=\B(\ell_2(\Lambda))\ox_{\min}C^*(\G)
\subset\B(\ell_2(\Lambda)\otimes\hh).$$ 
Composing these three morphisms, we obtain the conclusion. 
\end{proof} 

We fix a unitary $\upsilon_p$ in $\ell_\infty(\Lambda_p)$ for each $p$ 
with $(\upsilon_p\chi_p\mid\chi_p)=0$ 
and define an element $\upsilon_\omega$ in $\ell_\infty(\Lambda)$, 
for each subset $\omega$ of $\Omega$, by 
$$\upsilon_\omega(x)=\left\{
\begin{array}{ccc}
\upsilon_p(x) & \mbox{ if }x\in\Lp\mbox{ with }p\in\omega,\\
0 & \mbox{ if }x\in\Lp\mbox{ with }p\notin\omega.
\end{array}\right.$$
Let $\Sg=\{I,h,k\}$ be the finite set of generators of $\G$ 
and let $E_\omega$, for each subset $\omega$ of $\Omega$, 
be the four dimensional operator space in 
$\ell_\infty(\Lambda)\rtimes\G$ spanned by 
$\upsilon_\omega$ and $U(g)$, $g\in\Sg$. 

\begin{lem}\label{main}
The subset $\{E_\omega : \omega \mbox{ a subset of }\Omega\}$ 
of $\OS_4$ is non-separable in the cb Banach-Mazur distance topology.
\end{lem}
\begin{proof}
To prove $\{E_\omega\}_\omega$ is non-separable, 
suppose the contrary that $\{E_\omega\}_\omega$ is separable.
We fix $\e>0$. 
As it was argued in Remark 2.10 in \cite{jp}, 
it follows that one can find distinct $\omega$ and $\omega'$, 
and $q\in\omega'\setminus\omega$ such that 
there is a complete contraction $\p\colon E_\omega\to E_{\omega'}$ 
with $\| \upsilon_{\omega'}-\p(\upsilon_\omega)\|<\e$ 
and $\| U(g)-\p(U(g))\|<\e$ for $g\in\Sg$. 
Let $\tilde{\pi_q}\colon 
\ell_\infty(\Lambda)\rtimes\G\to\B(\ell_2(\Lq))$ be 
the covariant representation corresponding to the quotient map from 
$\ell_\infty(\Lambda)$ onto $\ell_\infty(\Lq)\subset\B(\ell_2(\Lq))$ 
and the unitary representation $\pi_q$ of $\G$ on $\ell_2(\Lq)$. 
It follows that $\psi=\tilde{\pi_q}\circ\p$ is a complete contraction 
with $\|\upsilon_q-\psi(\upsilon_\omega)\|<\e$ and 
$\|\pi_q(g)-\psi(U(g))\|<\e$ for $g\in\Sg$. 
Hence, we have that 
\begin{eqnarray*}
& & \hspace*{-1em}\|\overline{\upsilon_q}\otimes\upsilon_\omega+
 \sum_{g\in\Sg}\overline{\pi_q(g)}\otimes U(g)\|
 _{\B(\overline{\ell_2(\Lq)})\ox_{\min}(\ell_\infty(\Lambda)\rtimes\G)}\\
& & \geq\|\overline{\upsilon_q}\otimes\psi(\upsilon_\omega)
 +\sum_{g\in\Sg}\overline{\pi_q(g)}\otimes\psi(U(g))\|
 _{\B(\overline{\ell_2(\Lq)}\otimes\ell_2(\Lq))}\\
& & >\|\overline{\upsilon_q}\otimes\upsilon_q
 +\sum_{g\in\Sg}\overline{\pi_q(g)}\otimes\pi_q(g)\|
 _{\B(\overline{\ell_2(\Lq)}\otimes\ell_2(\Lq))}-4\e\\
& & =4(1-\e).
\end{eqnarray*}
Combining the above inequality with Lemma \ref{rep}, we have 
$$\|\overline{\upsilon_q}\otimes\upsilon_\omega\otimes 1+
 \sum_{g\in\Sg}\overline{\pi_q(g)}\otimes \pi(g)\otimes u(g)\|
 _{\B(\overline{\ell_2(\Lq)}\otimes\ell_2(\Lambda)\otimes\hh)}
> 4(1-\e)$$
By the uniform convexity of Hilbert spaces, we can find 
a unit vector 
$$\zeta=\sum_{(x,y)\in\Lq\times\Lambda}
\overline{\delta_x}\otimes\delta_y\otimes\zeta(x,y)
\in\overline{\ell_2(\Lq)}\otimes\ell_2(\Lambda)\otimes\hh$$ 
such that 
$$\|\zeta
-(\overline{\upsilon_q}\otimes\upsilon_\omega\otimes 1)\zeta\|<\beta(\e)
\mbox{ and }
\|\zeta
-(\overline{\pi_q(g)}\otimes \pi(g)\otimes u(g))\zeta\|<\beta(\e)$$ 
for $g\in\Sg$, where $\beta\colon\R_{\geq0}\to\R_{\geq0}$ 
is a continuous function with $\beta(0)=0$ 
(cf.\ the remarks following Theorem \ref{selberg}).
>From the first inequality and the fact that $q\notin\omega$, 
we may assume that $\zeta$ is zero on $\Lq\times\Lq$. 
We put 
$$\xi=\sum_{(x,y)\in\Lq\times\Lambda}
\overline{\delta_x}\otimes\delta_y\|\zeta(x,y)\|
\in\overline{\ell_2(\Lq)}\otimes\ell_2(\Lambda).$$ 
It follows that 
\begin{eqnarray*}
\hspace*{-1em}\beta(\e)
&>& \|\zeta
 -(\overline{\pi_q(g)}\otimes \pi(g)\otimes u(g))\zeta\|
 _{\overline{\ell_2(\Lq)}\otimes\ell_2(\Lambda)\otimes\hh}\\
&=& \|\hspace*{-0.7em}\sum_{(x,y)\in\Lq\times\Lambda}\hspace*{-0.7em}
 \overline{\delta_x}\otimes\delta_y\otimes
 \Bigl(\zeta(x,y)-u(g)\zeta(\sigma_q(g^{-1})x,\sigma(g^{-1})y)\Bigr)\|
 _{\overline{\ell_2(\Lq)}\otimes\ell_2(\Lambda)\otimes\hh}\\
&\geq& \|\hspace*{-0.7em}\sum_{(x,y)\in\Lq\times\Lambda}\hspace*{-0.7em}
 \overline{\delta_x}\otimes\delta_y\Bigl(\|\zeta(x,y)\|
 -\|u(g)\zeta(\sigma_q(g^{-1})x,\sigma(g^{-1})y)\|\Bigr)\|
 _{\overline{\ell_2(\Lq)}\otimes\ell_2(\Lambda)}\\
&=& \|\xi-(\overline{\pi_q(g)}\otimes \pi(g))\xi\|
 _{\overline{\ell_2(\Lq)}\otimes\ell_2(\Lambda)}
\end{eqnarray*}
for $g\in\Sg$. 
We are now in position to employ Selberg's theorem. 
Indeed, the unit vector $\xi$ is zero on $\Lq\times\Lq$ 
and, for $p\neq q$, Lemma \ref{irr} implies that 
any fixed vector for the representation $\overline{\pi_q}\otimes\pi_p$ 
is a scalar multiple of $\overline{\chi_q}\otimes\chi_p$.
Thus, it follows from Theorem \ref{selberg} that 
$$\|\xi-\overline{\chi_q}\otimes\eta
\|_{\overline{\ell_2(\Lq)}\otimes\ell_2(\Lambda)}
<\alpha(\beta(\e))$$
for some unit vector $\eta\in\ell_2(\Lambda)$, 
where $\alpha\colon\R_{\geq0}\to\R_{\geq_0}$ is a continuous function 
with $\alpha(0)=0$ 
(cf.\ the remarks following Theorem \ref{selberg}).
Finally, we have 
\begin{eqnarray*} 
\beta(\e)
&>& \|\zeta-(\overline{\upsilon_q}\otimes\upsilon_\omega\otimes 1)\zeta\|\\
&=& \|\xi-(\overline{\upsilon_q}\otimes\upsilon_\omega)\xi\|\\
&>& \|(\overline{\chi_q}\otimes\eta)
 -(\overline{\upsilon_q}\otimes\upsilon_\omega)
 (\overline{\chi_q}\otimes\eta)\|-2\alpha(\beta(\e))\\
&=& 2^{1/2}-2\alpha(\beta(\e)) 
\end{eqnarray*}
(recall that we chose $\upsilon_p$ so that 
$\chi_p\perp\upsilon_p\chi_p$), 
but this gives a contradiction when $\e>0$ is chosen sufficiently small. 
This completes the proof. 
\end{proof}

We now prove Proposition \ref{prop} and Theorem \ref{thm}. 

\begin{proof}
Proposition \ref{prop} follows from Lemma \ref{main} 
as we have explained in the first paragraph of this section. 
We turn to the proof of Theorem \ref{thm}.
To show that $\B(\ell_2)\ox_{\min}\B(\ell_2)$ 
does not have the WEP, suppose the contrary that 
$\B(\ell_2)\ox_{\min}\B(\ell_2)$ has the WEP. 
We take $\Lambda$ as in Proposition \ref{prop} so that, 
by Theorem \ref{kirchberg}, 
$$(\ell_\infty(\Lambda)\rtimes\G)\ox_{\min}\B(\ell_2)
\neq(\ell_\infty(\Lambda)\rtimes\G)\ox_{\max}\B(\ell_2).$$
By universality, we have the canonical isomorphism 
$$(\ell_\infty(\Lambda)\rtimes\G)\ox_{\max}\B(\ell_2)
=(\ell_\infty(\Lambda)\ox_{\min}\B(\ell_2))\rtimes\G,$$ 
here $\G$ acts on $\B(\ell_2)$ trivially. 
By the assumption that $\B(\ell_2(\Lambda))\ox_{\min}\B(\ell_2)$ 
has the WEP, arguing as Lemma \ref{rep}, 
one can show that the canonical morphism 
$$(\ell_\infty(\Lambda)\ox_{\min}\B(\ell_2))\rtimes\G
\subset\B(\ell_2(\Lambda)\otimes\ell_2\otimes\hh)$$ 
is faithful. 
Composing these, we have a faithful representation 
$$(\ell_\infty(\Lambda)\rtimes\G)\ox_{\max}\B(\ell_2)
\to\B(\ell_2(\Lambda)\otimes\hh\otimes\ell_2)$$ 
which obviously factors through 
$(\ell_\infty(\Lambda)\rtimes\G)\ox_{\min}\B(\ell_2)$. 
This is absurd. 
\end{proof} 
\appendix\section{A pathology in equivariant $KK$-theory}
With the help of Wassermann's construction \cite{wassermann}, 
we prove the following. 
\begin{thm2}\label{kex} 
Let $\G=SL(3,\Z)$ and let $M=\prod_{n=1}^\infty\M_n$. 
There is a short exact sequence of separable commutative 
$\G$-$C^*$-algebras $0\to J \to B \to A\to 0$
such that the corresponding sequence 
$$K_0(M\ox_{\min}(J\rtimes\G))
\to K_0(M\ox_{\min}(B\rtimes\G))
\to K_0(M\ox_{\min}(A\rtimes\G))$$ 
is not exact. 
\end{thm2}
\begin{cor2} 
The six-term exact sequence in $\G$-equivariant $KK$-theory 
fails to hold for the short exact sequence of separable commutative 
$\G$-$C^*$-algebras appearing in Theorem \ref{kex}. 
\end{cor2}
This corollary was pointed out by Skandalis and 
the proof is almost the same as \cite{skandalis}. 
Maghfoul \cite{maghfoul} proved that such a pathology does not occur 
under a certain $K$-theoretical amenability condition on $\G$. 
On the other hand, Higson, Lafforgue and Skandalis \cite{hls} found 
a similar pathology for Gromov's non-exact group $\G$ \cite{gromov}. 

 From now on, we denote the group $SL(3,\Z)$ by $\G$. 
For each prime number $p$, 
we let $\Lp$ be the projective space 
$(\Z_p^3-\{0\})/\Z_p^\times$ 
and observe that $|\Lp|=(p^3-1)/(p-1)=p^2+p+1$. 
The action of $\G$ on $\Z_p^3$ 
(through linear transformation by $SL(3,\Z_p)$) 
induces a transitive action $\sigma_p$ of $\G$ 
on the set $\Lp$. 
Let $\pi_p$ be the corresponding unitary representation 
of $\G$ on $\ell_2(\Lp)$. 
Let $\chi_p=|\Lp|^{-1/2}\sum_{x\in\Lp}\delta_x\in\ell_2(\Lp)$ 
be the constant function of norm $1$ and let $z\in\B(\ell_2(\Lp))$ 
be the projection onto the one dimensional subspace $\C\chi_p$. 
We observe that $z$ is in the center of $C^*(\pi_p(\G))$. 
Let $\pi_p^\circ$ be the subrepresentation of $\pi_p$, 
which is the restriction to the subspace $(1-z)\ell_2(\Lp)$. 
A similar proof to that of Lemma \ref{irr} in Section \ref{sec2} 
yields the following lemma. 
\begin{lem2}\label{irr3}
The representation $\pi_p^\circ$ is irreducible 
for every prime number $p$. 
\end{lem2} 
Let $\Omega$ be an infinite set of odd prime numbers 
such that $p,q\in\Omega$ and $p>q$ implies that $p\equiv 1\mod q$ 
(and in particular $p>2q$). 
Such an infinite set $\Omega$ exists by Dirichlet's theorem. 
For each $p\in\Omega$, we define the subset $X_p\subset\Lp$ by  
$$X_p=\{ (1 \ a \ b)^T\in\Lp : a=0,2,4,\ldots,p-1,\ b\in\Z_p\}$$
and observe that $|\Lp|/3<|X_p|=(p^2+p)/2<|\Lp|/2$. 
For $h=I_3+e_{21}\in\G$, we have that 
$$\sigma_p(h^q)X_p\cap X_p\cap\sigma_p(h^{-q})X_p=\emptyset$$ 
whenever $p,q\in\Omega$ are distinct.  
Indeed, this easily follows from the fact that 
we have either $q\equiv 1\mod p$ (and hence $\sigma_p(h^q)=\sigma_p(h)$) 
or $2q<p$ when $p,q\in\Omega$ are distinct. 
The action $\sigma_p$ induces an action of $\G$ on $\ell_\infty(\Lp)$ 
and on $\prod_{p\in\Omega}\ell_\infty(\Lp)$. 
We often identify $\ell_\infty(\Lp)$ with 
the diagonal of $\B(\ell_2(\Lp))$. 
For each $p\in\Omega$, let $e_p\in\ell_\infty(\Lp)$ be 
the characteristic function of $X_p$. 
Let $e=(e_p)_p\in\prod_{p\in\Omega}\ell_\infty(\Lp)$ be a projection, 
$B$ be the unital separable $\G$-$C^*$-subalgebra of 
$\prod_{p\in\Omega}\ell_\infty(\Lp)$ generated by 
$e$ and $J=\bigoplus_{p\in\Omega}\ell_\infty(\Lp)$ 
and let $Q$ be the $\G$-equivariant quotient from $B$ onto $A=B/J$. 
We denote by $U(g)$ the implementing unitary of $g\in\G$ in 
the full crossed product $B\rtimes\G$.

For $M=\prod_{q\in\Omega}\B(\overline{\ell_2(\Lq)})$ and 
a fixed finite symmetric set $\Sg$ of generators of $\G$ 
containing the unit, we define a selfadjoint element $s$ 
and a projection $t$ in $M\ox_{\min}(B\rtimes\G)$ by 
$$s=\frac{1}{|\Sg|}\sum_{g\in\Sg}\overline{\pi(g)}\otimes U(g) 
\mbox{ and }
t=\overline{e}\otimes e+(\overline{1-e})\otimes(1-e),$$
where $\pi(g)=(\pi_q(g))_q\in\prod_{q\in\Omega}\B(\ell_2(\Lq))$. 
Since $\G=SL(3,\Z)$ has the Kazhdan property \cite{hv} \cite{valetteT}, 
there is $0<\e<1$ such that $\spec(s)\subset [-1,1-\e]\cup\{1\}$
(cf.\ the remarks following Theorem \ref{selberg}). 
We will prove that the spectrum of $r=s+t$ has a gap around 2. 
For this reason, we decompose $r$ into a direct sum. 
For each $q\in\Omega$, let $Q_q$ be 
the canonical quotient from $\prod_{p\in\Omega}\ell_\infty(\Lp)$ 
onto $\ell_\infty(\Lq)$, and let $Q_q'$ be the canonical quotient 
from $\prod_{p\in\Omega}\ell_\infty(\Lp)$ 
onto $\prod_{p\in\Omega\setminus\{q\}}\ell_\infty(\Lq)$. 
We still denote their restriction to $B$ by $Q_p$ and $Q_q'$ and 
let $A_q=Q_p(B)=\ell_\infty(\Lq)$ and $A_q'=Q_p'(B)$. 
Since we have $B=A_q\oplus A_q'$ as a $\G$-$C^*$-algebra, 
$B\rtimes\G=A_q\rtimes\G\oplus A_q'\rtimes\G$. 
Therefore, we have 
$$M\ox_{\min}(B\rtimes\G)\subset
\prod_{q\in\Omega}\left(
\B(\overline{\ell_2(\Lq}))\ox_{\min}(A_q\rtimes\G)
\oplus\B(\overline{\ell_2(\Lq}))\ox_{\min}(A_q'\rtimes\G)\right).$$
We denote respectively by $r_q$, $s_q$ and $t_q$ the direct summands of 
$r$, $s$ and $t$ in 
$\B(\overline{\ell_2(\Lq}))\ox_{\min}(A_q\rtimes\G)$
and by $r_q'$, $s_q'$ and $t_q'$ the direct summands of 
$r$, $s$ and $t$ in 
$\B(\overline{\ell_2(\Lq}))\ox_{\min}(A_q'\rtimes\G)$. 
\begin{lem2}\label{gap}
We have $2\in\spec(r_q)\subset[-1,2-10^{-4}\e]\cup\{2\}$ and 
$\spec(r_q')\subset[-1,2-10^{-4}\e]$. 
\end{lem2}
\begin{proof} 
Let $C^*(\G)\subset\B(\hh)$ be a faithful representation and 
denote by $u(g)$ the unitary in $C^*(\G)$ corresponding to $g\in\G$. 
Since $A_q=\ell_\infty(\Lq)$ is finite dimensional, 
it is not hard to see that the representation 
of $A_q\rtimes\G$ on $\ell_2(\Lq)\otimes\hh$ 
given by $A_q\ni a\mapsto a\otimes 1$ and 
$\G\ni g\mapsto \pi_q(g)\otimes u(g)$ 
is faithful. 
Hence, we have 
$$s_q=\frac{1}{|\Sg|}\sum_{g\in\Sg}
\overline{\pi_q(g)}\otimes\pi_q(g)\otimes u(g) 
\mbox{ and }
t_q=\overline{e_q}\otimes e_q\otimes 1
+(\overline{1-e_q})\otimes(1-e_q)\otimes 1$$
on $\overline{\ell_2(\Lq)}\otimes\ell_2(\Lq)\otimes\hh$. 
We identify the Hilbert space 
$\overline{\ell_2(\Lq)}\otimes\ell_2(\Lq)\otimes\hh$
with the space $\hk$ of Hilbert-Schmidt operators 
from $\ell_2(\Lq)$ to $\ell_2(\Lq)\otimes\hh$ 
so that $\overline{\pi_q(g)}\otimes\pi_q(g)\otimes u(g)$ 
acts on $\hk$ by 
$\hk\ni T\mapsto (\pi_q(g)\otimes u(g))T\pi_q(g)^*\in\hk$.  
Then, it follows from the uniform convexity of a Hilbert space 
(cf. the remarks following Theorem \ref{selberg}) 
that the eigenspace $\hk_0$ of $s_q$ w.r.t.\ the eigenvalue $1$ is 
$$\hk_0=\{ T\in\hk : 
(\pi_q(g)\otimes u(g))T\pi_q(g)^*=T\mbox{ for all }g\in\G\}.$$
For each $\xi\colon\Lq\to\hh$, 
we associate $S_\xi\in\hk$ 
defined by $S_\xi(\delta_x)=\delta_x\otimes\xi_x$ 
and then we define the subspace $\hk_1$ of $\hk$ by 
$$\hk_1=\{ S_\xi\in\hk : 
\xi\mbox{ satisfies }u(g)\xi_x=\xi_{\sigma_q(g)x}
\mbox{ for all }x\in\Lq\mbox{ and }g\in\G\}.$$
Since $\hh$ contains a non-zero fixed vector for 
the unitary representation $u$, it is not too difficult to see that 
$\hk_1$ is non-empty and contained 
in the intersection of eigenspaces of $s_q$ and $t_q$ 
w.r.t.\ the eigenvalues $1$ (and hence $r_q|_{\hk_1}=2$). 
We claim that $\spec(r_q|_{\hk_1^\perp})\subset[-1,2-10^{-4}\e]$. 
This easily follows if we prove that $\| t_q(T)\|<(25/36)^{1/2}\| T\|$ 
for any $T\in\hk_0\ominus\hk_1$. 
To prove this, we give ourselves $T\in\hk_0\ominus\hk_1$ of norm 1. 
Since $T\in\hk_0$, we have that $\xi\colon\Lq\to\hh$, given by 
$\xi_x=(T\delta_x)(x)$, satisfies $u(g)\xi_x=\xi_{\sigma_q(g)x}$ 
for every $x\in\Lq$ and $g\in\G$. 
It follows that $0=(T,S_\xi)_{\hk}=\sum_x\|(T\delta_x)(x)\|^2$ 
and hence $(T\delta_x)(x)=0$ for all $x\in\Lq$. 
We define $\tilde{T}\in\B(\ell_2(\Lq))$ by 
$\tilde{T}\delta_x=|T\delta_x|$, 
where for $\zeta\in\ell_2(\Lq)\otimes\hh$, 
the vector $|\zeta|\in\ell_2(\Lq)$ is given by 
$|\zeta|(x)=\|\zeta(x)\|_{\hh}$. 
Since $\tilde{T}$ commutes with $\pi_q(\G)$ and 
its diagonal entries are zero, it follows from Lemma \ref{irr3} that 
$\tilde{T}=\lambda1+\mu z$ with $n\lambda+\mu=0$, where $n=|\Lq|$. 
Since $\|\tilde{T}\|=\| T\|=1$, we obtain 
$\tilde{T}=(n(n-1))^{-1/2}(nz-1)$. 
Therefore, we have 
\begin{eqnarray*} 
\| t_q(T)\|^2 
&=& \| (e_q\otimes 1)Te_q+((1-e_q)\otimes 1)T(1-e_q)\|_{\hk}^2\\
&=& \sum_{x\in X_q}\| e_q\tilde{T}(\delta_x)\|^2 
 + \sum_{x\in\Lq\setminus X_q}\|(1-e_q)\tilde{T}(\delta_x)\|^2\\
&=& (n(n-1))^{-1}(|X_q|(|X_q|-1)+(n-|X_q|)(n-|X_q|-1))\\
&<& (1/2)^2+(2/3)^2=25/36
\end{eqnarray*}
since $n/3<|X_q|<n/2$. 
This completes the proof of the first half. 

Take a faithful representation $A_q'\rtimes\G\subset\B(\hh)$ and  
denote $e_q'=Q_q'(e)\in A_q'$ and 
$U_q'(g)\in A_q'\rtimes\G$ the implementing unitary for $g\in\G$. 
Then, by the construction, we have 
$$f:=3^{-1}(\Ad U_q'(h^q)(e_q')+e_q'+\Ad U_q'(h^{-q})(e_q'))\le 2/3$$ 
in $A_q'$ as 
$\sigma_p(h^q)X_p\cap X_p\cap\sigma_p(h^{-q})X_p=\emptyset$ 
for all $p\in\Omega\setminus\{q\}$. 
We note that 
$$s_q'=\frac{1}{|\Sg|}\sum_{g\in\Sg}\overline{\pi_q(g)}\otimes U_q'(g) 
\mbox{ and }
t_q'=\overline{e_q}\otimes e_q'+(\overline{1-e_q})\otimes(1-e_q')$$
on $\overline{\ell_2(\Lq)}\otimes\hh$. 
We identify the Hilbert space $\overline{\ell_2(\Lq)}\otimes\hh$ 
with the space $\hk$ of Hilbert-Schmidt operators 
from $\ell_2(\Lq)$ to $\hh$ 
so that $\overline{\pi_q(g)}\otimes U_q'(g)$ 
acts on $\hk$ by $\hk\ni T\mapsto U_q'(g)T\pi_q(g)^*\in\hk$.  
Then, it follows from the uniform convexity of a Hilbert space 
(cf. the remarks following Theorem \ref{selberg}) 
that the eigenspace $\hk_0$ of $s_q'$ w.r.t.\ the eigenvalue $1$ is 
$$\hk_0=\{ T\in\hk : U_q'(g)T\pi_q(g)^*=T\mbox{ for all }g\in\G\}.$$
We claim that $\| t_q'(T)\|<(8/9)^{1/2}\| T\|$ for any $T\in\hk_0$. 
Then the second half of this lemma follows. 
To prove this, we give ourselves $T\in\hk_0$ of norm $1$.
Since $U_q'(h^q)T=T\pi_q(h^q)=T$, we have 
\begin{eqnarray*} 
\| t_q(T)\|^2 
&=& \Tr(T^*e_q'Te_q+T^*(1-e_q')T(1-e_q))\\ 
&=& \Tr(T^*fTe_q+T^*(1-f)T(1-e_q))\\ 
&\le& 1-3^{-1}\Tr(T^*Te_q).
\end{eqnarray*}
Since $T^*T$ commutes with $\pi_q(\G)$, 
it follows from Lemma \ref{irr3} that 
$T^*T=\lambda1+\mu z$ for some real number $\lambda$ and $\mu$ 
with $\lambda|\Lq|+\mu=1$.  
Hence, we have that 
$\Tr(T^*Te_q)=|X_q|/|\Lq|>1/3$. 
This completes the proof.
\end{proof} 

We now prove Theorem \ref{kex}. 

\begin{proof}
Since we have $\spec(r)\subset[-1,2-10^{-4}\e]\cup\{2\}$ 
by Lemma \ref{gap}, the spectral projection $d$ of $r$ 
corresponding to the spectral subset $\{2\}$ is contained 
in the $C^*$-algebra $M\ox_{\min}(B\rtimes\G)$. 
Since 
$$M\ox_{\min}(A\rtimes\G)\subset\prod_{q\in\Omega}\left(
\B(\overline{\ell_2(\Lq}))\ox_{\min}(A\rtimes\G)\right)$$
and the quotient $Q$ from $B$ onto $A$ factors through $A_q'$ 
for each $q\in\Omega$, 
we have that $(\id_M\otimes\tilde{Q})(d)=0$ by Lemma \ref{gap}, 
where $\tilde{Q}$ is the quotient 
from $B\rtimes\G$ onto $A\rtimes\G$ induced by $Q$. 
Finally, we observe that 
the $K_0$-element corresponding to $d$ does not come from 
$K_0(M\ox_{\min}(J\rtimes\G))$ as 
any element from $K_0(M\ox_{\min}(J\rtimes\G))$ vanishes 
on $\tau_q$ for all but finitely many $q\in\Omega$, 
where $\tau_q$ is the tracial state on 
$\B(\overline{\ell_2(\Lq)})\ox_{\min}\B(\ell_2(\Lq))$
evaluated through $\tilde{\pi_q}\colon B\rtimes\G\to\B(\ell_2(\Lq))$. 
\end{proof}

\end{document}